\documentclass[10pt,letterpaper]{amsart}

\usepackage[margin=1.5in]{geometry}

\usepackage{standalone}

\usepackage{amsrefs}
\usepackage{amsfonts}
\usepackage{amssymb}
\usepackage{amsmath}
\usepackage{amsthm}
\usepackage[inline]{enumitem}
\usepackage{hyperref}
\usepackage{enumitem}

\usepackage{color}

\usepackage[textsize=small]{todonotes}
\newcommand{\N}{\mathbb{N}}
\newcommand{\R}{\mathbb{R}}
\newcommand{\C}{\mathbb{C}}
\newcommand{\st}{ : }

\renewcommand{\>}{\rangle}
\newcommand{\restrict}{\upharpoonright}

\newtheorem{theorem}{Theorem}

\newtheorem{lemma}[theorem]{Lemma}

\newtheorem{question}[theorem]{Question}

\title{Sparse analytic systems}

\author{Brent Cody}
\address{
Department of Mathematics and Applied Mathematics \\
Virginia Commonwealth University \\
1015 Floyd Avenue \\
Richmond, Virginia 23284, USA}
\email{bmcody@vcu.edu}

\author{Sean Cox}
\address{
Department of Mathematics and Applied Mathematics \\
Virginia Commonwealth University \\
1015 Floyd Avenue \\
Richmond, Virginia 23284, USA.}
\email{scox9@vcu.edu}

\author{Kayla Lee}
\address{
Department of Mathematics and Applied Mathematics \\
Virginia Commonwealth University \\
1015 Floyd Avenue \\
Richmond, Virginia 23284, USA.  }
\email{leek10@vcu.edu}

\thanks{We gratefully acknowledge suggestions of the two referees.  Partially funded by VCU SEED Grant and NSF grant DMS-2154141 of the second author.}

\subjclass[2020]{03E50,  03E25, 26E05, 30D20 }

\keywords{Wetzel's Problem, analytic functions, Continuum Hypothesis}

\begin{document}

\begin{abstract}
Erd\H{o}s~\cite{MR168482} proved that the Continuum Hypothesis (CH) is equivalent to the existence of an uncountable family $\mathcal{F}$ of (real or complex) analytic functions, such that $\big\{ f(x) \ : \ f \in \mathcal{F} \big\}$ is countable for every $x$.  We strengthen Erd\H{o}s' result by proving that CH is equivalent to the existence of what we call \emph{sparse analytic systems} of functions.  We use such systems to construct, assuming CH, an equivalence relation $\sim$ on $\mathbb{R}$ such that any ``analytic-anonymous" attempt to predict the map $x \mapsto [x]_\sim$ must fail almost everywhere.  This provides a consistently negative answer to a question of Bajpai-Velleman~\cite{MR3552748}.  
\end{abstract}

\maketitle

\section{Introduction}\label{sec_Intro}

In the early 1960s, John Wetzel posed:
\begin{quote}
\textbf{Wetzel's Problem:}  
If $\mathcal{F}$ is  a family of analytic functions (on some common domain) such that $\{ f(x) \ : \ f \in \mathcal{F} \}$ is countable for every $x$, must $\mathcal{F}$ be a countable family?  
\end{quote}

A few years later, Erd\H{o}s proved that an affirmative answer to Wetzel's Problem is equivalent to the negation of Cantor's \emph{Continuum Hypothesis (CH)}.  Combined with Paul Cohen's proof of the independence of CH, this showed that Wetzel's Problem is independent of the standard axioms of mathematics (ZFC).  Upon learning of Erd\H{o}s' theorem, Wetzel remarked to his dissertation advisor (Halsey Royden) that ``\dots once again a
natural analysis question has grown horns!"  This quote, and other interesting history surrounding Wetzel's Problem, appears in Garcia-Shoemaker~\cite{GarciaShoemaker}.  Erd\H{o}s' proof even made it into Aigner-Ziegler's ``Proofs from the Book" (\cite{MR3823190}).  It will be more convenient for us to state and refer to Erd\H{o}s' equivalence in the negated form:
\begin{theorem}[Erd\H{o}s~\cite{MR168482}]\label{thm_Erdos}
The following are equivalent:
\begin{enumerate}
    \item CH;
    \item\label{item_WetzelSystem} There exists an uncountable family $\mathcal{F}$ of analytic functions on some fixed open domain $D$ of either $\mathbb{R}$ or $\mathbb{C}$, such that for every $x \in D$, 
    \[
    \big\{  f(x) \ : \ f \in \mathcal{F} \big\}
    \]
    is countable.
\end{enumerate}
\end{theorem}

Motivated by connections to work of Hardin-Taylor (\cite{MR2384262}, \cite{MR3100500}) and Bajpai-Velleman~\cite{MR3552748} described below, we strengthen Theorem \ref{thm_Erdos} as follows.  If $P \in \mathbb{R}^2$ we denote the first coordinate of $P$ by $x_P$ and the second coordinate by $y_P$.  Define a \textbf{sparse (real) analytic system} to mean a collection 
\[
\big\{ f_P \ : \ P  \in \mathbb{R}^2  \big\}
\]
such that:
\begin{enumerate}
    \item for all $P \in \mathbb{R}^2$, $f_P$ is an increasing, analytic bijection from $\mathbb{R} \to \mathbb{R}$ that passes through the point $P$; and

    \item For all $z \in \mathbb{R}$, the sets
    \[
    \big\{ f_P(z) \ : \ P \in \mathbb{R}^2 \text{ and }  z \ne x_P \big\}
    \]
    and
    \[
    \big\{ f^{-1}_P(z) \ : \ P \in \mathbb{R}^2 \text{ and } z \ne y_P \big\}
    \]
    are both countable.
\end{enumerate}
We prove the following strengthening of Erd\H{o}s' Theorem \ref{thm_Erdos}:
\begin{theorem}\label{thm_StrengthenErdos}
The following are equivalent:
\begin{enumerate}
    \item CH
    \item There exists a sparse real analytic system.

\end{enumerate}
\end{theorem}

We use Theorem \ref{thm_StrengthenErdos} to answer a question of Bajpai and Velleman, assuming CH.  Given a nonempty set $S$, let $\boldsymbol{{}^{\mathbb{R}} S}$ denote the collection of total functions from $\mathbb{R}$ to $S$, and $\boldsymbol{{}^{\underset{\smile}{\mathbb{R}}} S}$ denote the collection of all $S$-valued functions $f$ such that $\text{dom}(f)=(-\infty,t_f)$ for some $t_f \in \mathbb{R}$.  An \textbf{$\boldsymbol{S}$-predictor} will refer to any function $\mathcal{P}$ with domain and codomain as follows:
\begin{equation}\label{eq_DomCodP}
\mathcal{P}: {}^{\underset{\smile}{\mathbb{R}}} S \to S.
\end{equation}
An $S$-predictor $\mathcal{P}$ will be called \textbf{good} if for all $F \in {}^{\mathbb{R}} S$, the set
\[
\Big\{ t \in \mathbb{R} \ : \  F(t) = \mathcal{P} \big( F \restriction (-\infty,t)  \big)  \Big\}
\]
has full measure in $\mathbb{R}$.  So $\mathcal{P}$ is good if for any total $F:\mathbb{R} \to S$, $\mathcal{P}$ ``almost always" correctly predicts $F(t)$ based only on $F \restriction (-\infty,t)$.\footnote{Note that $F$ is allowed to be highly discontinuous; otherwise the problem trivializes, since one could simply predict $F(t)$ by considering $\lim_{x \nearrow t} F(x)$, which only depends on $F \restriction (-\infty,t)$.}   Hardin-Taylor~\cite{MR2384262} proved that for any set $S$, there exists a good $S$-predictor, and in \cite{MR3100500} they raised the question of whether these good predictors could also be arranged to be ``$\Gamma$-anonymous" with respect to certain classes $\Gamma \subseteq \text{Homeo}^+(\mathbb{R})$;\footnote{$\text{Homeo}^+(\mathbb{R})$ denotes the set of increasing homeomorphisms from $\mathbb{R}$ to $\mathbb{R}$. } an $S$-predictor $\mathcal{P}$ is $\Gamma$-anonymous if for every $\varphi \in \Gamma$ and every $f \in {}^{\underset{\smile}{\mathbb{R}}} S$,
\[
\mathcal{P}\big( f \big) = \mathcal{P} \big( f \circ \varphi \big), 
\]
where $f \circ \varphi$ is the member of ${}^{\underset{\smile}{\mathbb{R}}} S$ whose domain is understood to be $\big(-\infty, \varphi^{-1}(t_f) \big)$.  Bajpai and Velleman~\cite{MR3552748} gave a positive and a negative result:
\begin{itemize}
    \item For every set $S$, there exists a good $S$-predictor that is anonymous with respect to the class of affine functions on the reals.  This strengthened a previous theorem of Hardin-Taylor~\cite{MR3100500}, who had gotten the same result for the smaller class of affine functions of slope 1 (i.e, shifts).
    \item There is an equivalence relation $\sim$ on $\mathbb{R}$ such that, letting $S:=\mathbb{R}/\sim$, there is \textbf{no} good $S$-predictor that is anonymous with respect to the class of increasing $C^\infty$ bijections on $\mathbb{R}$.
\end{itemize}
They asked about classes intermediate between the affine functions and the $C^\infty$ functions; in particular:
\begin{question}[Bajpai-Velleman~\cite{MR3552748}, page 788]\label{q_BV}
Does there exist (for every set $S$) a good $S$-predictor that is anonymous with respect to the analytic members of $\text{Homeo}^+(\mathbb{R})$?
\end{question}

We use Theorem \ref{thm_StrengthenErdos}, together with an argument from Bajpai-Velleman~\cite{MR3552748}, to prove:
\begin{theorem}\label{thm_AnswerBV_CH}
Assuming CH, the answer to Question \ref{q_BV} is negative.
\end{theorem}

Section \ref{sec_Franklin} provides an interpolation theorem that will be used in the proof of Theorem \ref{thm_StrengthenErdos},  Section \ref{sec_Erdos} proves Theorem \ref{thm_StrengthenErdos}, Section \ref{sec_BV_solution} proves Theorem \ref{thm_AnswerBV_CH}, and Section \ref{sec_Concluding} has concluding remarks and open questions.

\section{An interpolation theorem}\label{sec_Franklin}

A key part of the proof of Theorem \ref{thm_StrengthenErdos} is the (ZFC) Theorem \ref{thm_StrengthenFranklin} below.  One of the referees pointed out that Theorem \ref{thm_StrengthenFranklin} follows from known results; in particular, it follows from the much more powerful Theorem 3.2 of Burke~\cite{MR3563069} or,  with modifications in the proofs, either Theorem 2 of Barth-Schneider~\cite{MR269834} or Corollary 1.9 of Burke~\cite{MR2485411}.  Since deriving Theorem \ref{thm_StrengthenFranklin} from those more powerful theorems is not trivial, we choose to present our original direct proof of Theorem \ref{thm_StrengthenFranklin}.

Recall that Cantor proved that any two countable dense subsets of $\mathbb{R}$ are order-isomorphic; and that this order-isomorphism easily extends uniquely to a homeomorphism of $\mathbb{R}$.  Franklin~\cite{MR1501300} considered the question of how nice this homeomorphism could be arranged to be, and showed that if $D$ and $E$ are countable dense subsets of $\mathbb{R}$, there is an order-isomorphism of $D$ with $E$ that extends to a real analytic function.  A series of papers improved this result, culminating in Barth-Schneider~\cite{MR269834}, who proved that there is an order-isomorphism of $D$ with $E$ that extends to an entire function $f: \mathbb{C} \to \mathbb{C}$, answering (one interpretation of) Question 24 of Erd\H{o}s~\cite{MR98702}.\footnote{See also Maurer~\cite{MR215994}, Nienhuys-Thiemann~\cite{MR0460638}, and Sato-Rankin~\cite{MR346157} for related results.  Burke~\cite{MR2485411} provides a nice historical overview of this literature on this topic.}  Subsequent work of Burke, mentioned above, further strengthened those results.  The variant we'll need for the proof of Theorem \ref{thm_StrengthenErdos} is:

\begin{theorem}\label{thm_StrengthenFranklin}
Suppose $\mathcal{D}$ is a partition of $\mathbb{R}$ into dense subsets of $\mathbb{R}$; for each $z \in \mathbb{R}$ let $D_z$ denote the unique $D \in \mathcal{D}$ such that $z \in D$.

Then for any $P=(x_P,y_P) \in \mathbb{R}^2$ and any countable set $W$ of reals, there is an entire function $f: \mathbb{C} \to \mathbb{C}$ such that:
\begin{enumerate}
    \item $f \restriction \mathbb{R}$ is real-valued (hence analytic, since $f:\mathbb{C} \to \mathbb{C}$ is entire), 
    \item $f \restriction \mathbb{R}$ is a bijection with strictly positive derivative;
    \item $f(x_P)=y_P$; and
    \item for each $w \in W$:
    \begin{enumerate}
        \item  if $w \ne x_P$ then $f(w) \in D_w$;
        \item  if $w \ne y_P$ then $f^{-1}(w) \in D_w$.
    \end{enumerate}
\end{enumerate}
\end{theorem}

Let us give a brief outline of the following proof of Theorem \ref{thm_StrengthenFranklin}, which is inspired by the proof of Nienhuys-Thiemann~\cite{MR0460638}. We will inductively define a sequence of functions $\<f_n\st n\in\N\>$ whose limit will be the desired function $f$. Each function $f_n$ will satisfy a version of Theorem \ref{thm_StrengthenFranklin}(4) for finitely many points in $W$. When we define the next function $f_{n+1}$ we will want it to be equal to $f_n$ on these finitely many points in $W$ that have already been taken care of, and we will want $f_{n+1}$ to satisfy Theorem \ref{thm_StrengthenFranklin}(4a) or Theorem \ref{thm_StrengthenFranklin}(4b), depending on whether $n$ is even or odd, for an additional point in $W$. We will write $A_n$ to denote the set of finitely many points of $W$ that have already been taken care of at stage $n$ with regard to Theorem \ref{thm_StrengthenFranklin}(4a), and we will write $B_n$ to denote the set of finitely many points of $W$ that have been taken care of in regard to Theorem \ref{thm_StrengthenFranklin}(4b).

Suppose $\mathcal{D}$ is a partition of $\mathbb{R}$ into dense sets, $W$ is a countable set of real numbers, and $P=(x_P,y_P)$ is a point in $\mathbb{R}^2$.  Fix a 1-1 enumeration $\{ w_n \ : \ n \in \mathbb{N} \}$ of $W$, and for each $n$ let $D_n$ be the unique member of $\mathcal{D}$ containing $w_n$.  Since $\mathcal{D}$ is a partition, we have
\begin{equation}\label{eq_KeyAssumption}
    \forall k,n \in \mathbb{N} \ \ \big( w_k \in D_n \ \iff \ D_k=D_n \ \iff \ w_n \in D_k \big).\tag{$*$}
\end{equation}
Suppose $p:\R\to\R$ is a continuous positive function such that
\begin{equation}\label{eq_GrowFast}
\forall n \in \mathbb{N} \     \lim_{t\rightarrow\infty}\frac{p(t)}{t^n}=\infty.
\end{equation}
We will inductively define sequences $\<f_n\st n\in\N\>$, $\<A_n\st n\in\N\>$ and  $\<B_n\st n\in\N\>$ such that $A_0=\emptyset$ and $B_0=\emptyset$ and for all $n\in\N$ we have
\begin{enumerate}
\item[(I)$_n$] $f_n:\C\to\C$ is entire and $f_n\restrict \R$ is real-valued;
\item[(II)$_n$] $f_n(x_P)=y_P$;
\item[(III)$_n$] $\forall x\in \R$ $f_n'(x)\geq\frac{1}{2}+\frac{1}{2^n}$ and thus $f_n\restrict \R$ is a bijection;
\item[(IV)$_n$] if $n>0$ then $\forall z\in\C$ $|f_n(z)-f_{n-1}(z)|<\frac{1}{2^n}p(|z|)$;  
\item[(V)$_n$] if $n=2k+1$ is odd then $A_n=A_{n-1}\cup\{w_k\}$, $B_n=B_{n-1}$ and we have $w_k\neq x_P$ $\implies$ $f_n(w_k)\in D_k$;
\item[(VI)$_n$] if $n=2k+2$ is even then $A_n=A_{n-1}$, $B_n=B_{n-1}\cup \{w_k\}$ and we have $w_k\neq y_P$ $\implies$ $f_n^{-1}(w_k)\in D_k$; and 
\item[(VII)$_n$] if $n>0$ then $f_n\restrict A_{n-1}=f_{n-1}\restrict A_{n-1}$ and $f_n^{-1}\restrict B_{n-1}=f_{n-1}^{-1}\restrict B_{n-1}$.
\end{enumerate}

First let us show that, assuming we have sequences $\<f_n\st n\in\N\>$ and $\<A_n\st n\in\N\>$ and $\<B_n\st n\in\N\>$ satisfying (I)$_n$-(VII)$_n$ for all $n$, the pointwise limit defined by $f(z)=\lim_{n\to\infty}f_n(z)$ has all of the desired properties. Suppose $D$ is any compact subset of $\C$. Since $\sum_{n=1}^\infty\frac{1}{2^n}$ converges and since $p(|z|)$ is bounded on $D$, the fact that (IV)$_n$ holds for all $n$ ensures that the sequence $\<f_n\st n\in\N\>$ is uniformly Cauchy on $D$. Hence we can define a function $f:\C\to\C$ by letting $f(z)=\lim_{n\to\infty}f_n(z)$. Since the sequence $\<f_n\st n\in\N\>$ is uniformly Cauchy on any compact set, it follows that the convergence of $\<f_n\st n\in\N\>$ to $f$ is uniform on any compact set, and hence $f$ is an entire function.

Now let us verify that Theorem \ref{thm_StrengthenFranklin}(1)-(4) hold for $f$. By (I)$_n$ and closure of $\R$ in $\C$ we see that $f\restrict \R$ is real valued and since (III)$_n$ holds for all $n$, we have $f'(x)\geq\frac{1}{2}$ for all $x\in\R$, thus Theorem \ref{thm_StrengthenFranklin}(1) and Theorem \ref{thm_StrengthenFranklin}(2) hold. Theorem \ref{thm_StrengthenFranklin}(3) holds since the sequence $\<f_n(x_P)\st n\in\N\>$ is constantly equal to $y_P$. To show that Theorem \ref{thm_StrengthenFranklin}(4) holds, let us prove that for all $i\in \N$ if $w_i\neq x_P$ then $f(w_i)\in D_i$ and if $w_i\neq y_P$ then $f^{-1}(w_i)\in D_i$. Fix $i\in\N$. We have $w_i\in A_{2i+1}$ and $w_i\in B_{2i+2}$ and furthermore, by (V)$_{2i+1}$ and (VI)$_{2i+2}$, $w_i\neq x_P$ implies $f_{2i+1}(w_i)\in D_i$ and $w_i\neq y_P$ implies $f_{2i+2}^{-1}(w_i)\in D_i$. Since (VII)$_n$ holds for all $n$ we see that both of the sequences $\<f_n(w_i)\st n\in\N\>$ and $\<f_n^{-1}(w_i)\st n\in\N\>$ are eventually constant, and indeed, for $n\geq 2i+2$ we have $f_n(w_i)=f_{2i+1}(w_i)$ and $f_n^{-1}(w_i)=f_{2i+2}^{-1}(w_i)$. Therefore $f(w_i)=f_{2i+1}(w_i)$ and $f^{-1}(w_i)=f_{2i+2}^{-1}(w_i)$, so (4) holds.

It remains to show that we can inductively define sequences $\<f_n\st n\in\N\>$, $\<A_n\st n\in\N\>$, and $\<B_n\st n\in\N\>$ that satisfy (I)$_n$-(VII)$_n$ for all $n\in\N$.

Let $f_0:\C\to\C$ be $f_0(z)=\frac{3}{2}(z-x_P)+y_P$, $A_0=\emptyset$ and $B_0=\emptyset$. One may easily verify that (I)$_0$-(VII)$_0$ hold.  For $n > 0$, Section \ref{sec_n_odd} shows how $f_n$ is constructed when $n$ is odd, and Section \ref{sec_n_even} shows how $f_n$ is constructed when $n$ is even.

\subsection{When $n$ is odd}\label{sec_n_odd}  Suppose $n=2k+1>0$ is odd and that $f_i$, $A_i$ and $B_i$ satisfying (I)$_i$-(VII)$_i$ have already been defined for $i\leq 2k$. If $k=0$ we have $A_0=\emptyset$ and $B_0=\emptyset$, whereas if $k>0$ we have
\[A_{n-1}=A_{2k}=A_{2(k-1)+2}=\{w_0,\ldots,w_{k-1}\}\] 
and
\[B_{n-1}=B_{2k}=\{w_0,\ldots,w_{k-1}\}.\]
In any case, we let $A_{n}=A_{n-1}\cup\{w_k\}$ and $B_{n}=B_{n-1}$. We define $f_{n} = f_{2k+1}$ in two cases as follows. 

\noindent \textbf{Case \ref{sec_n_odd}.A:  $\boldsymbol{w_k\notin \{x_P\}\cup A_{n-1}\cup f_{n-1}^{-1}(B_{n-1})}$}.   Let us argue that there is an entire function $g_{n}$ such that
\begin{enumerate}
\item[(i)] $(\forall z\in \C)$ $g_{n}(z)=0$ $\iff$ $z\in \{x_P\}\cup A_{n-1}\cup f_{n-1}^{-1}(B_{n-1})$,
\item[(ii)] $(\forall z\in \C)$ $|g_{n}(z)|\leq\frac{1}{2^{n}}p(|z|)$ and
\item[(iii)] $(\forall x\in \R)$ $g_{n}'(x)\geq -\frac{1}{2^{n}}$.
\end{enumerate}
Take
\[h_n(z)=(z-x_P)^{\beta_n}(z-w_0)\cdots(z-w_{k-1})(z-f_{n-1}^{-1}(w_0))\cdots(z-f_{n-1}^{-1}(w_{k-1})),\]
where $\beta_n\in\{1,2\}$ is such that the degree of $h_n$ is odd. We will show that for small enough positive $\alpha_n\in \R$ the function $g_n(z)=\alpha_nh_n(z)$ satisfies (i)-(iii). Clearly $h_n$ satisfies (i), so any such function $g_n(z)$ satisfies (i). For (ii), choose $m\in\N$ and some positive $c\in\R$ such that $|h_n(z)|\leq |z|^m+c$ for all $z\in \C$. By our assumption on $p$ we have $\lim_{|z|\rightarrow\infty}\frac{p(|z|)}{|z|^m+c}=\infty$, and thus we can let $D\subseteq\C$ be a large enough closed disk centered at the origin such that $z\in \C\setminus D$ implies $1\leq \frac{p(|z|)}{|z|^m+c}$. Since $p$ is a continuous positive function, we can choose a positive $\alpha_n\in\R$ such that $\alpha_n\leq\frac{1}{2^n}$ and $\alpha_n\leq\frac{p(|z|)}{2^n(|z|^m+c)}$ for all $z\in D$. Then it follows that for every $z\in\C$ we have
\[|\alpha_nh_n(z)|\leq\alpha_n(|z|^m+c)\leq\frac{1}{2^n}p(|z|).\]
Let us verify that (iii) holds for small enough $\alpha_n$. Since $h_n$ is odd and has a positive leading coefficient, the derivative of $h_n\restrict\R$ is bounded below. So we may let $d=\inf\{h_n'(x)\st x\in\R\}\in\R$. Thus we may choose a small enough positive $\alpha_n\in\R$ such that $\alpha_nd\geq-\frac{1}{2^n}$, and then it follows that for all $x\in \R$ we have $\alpha_nh_n'(x)\geq\alpha_n d\geq-\frac{1}{2^n}$.

Using the case assumption that $w_k\notin \{x_P\}\cup A_{n-1}\cup f_{n-1}^{-1}(B_{n-1})$, we see that $g_n(w_k)\neq 0$ and hence it follows that the set
\[\{f_{n-1}(w_k)+Mg_{n}(w_k)\st M\in[0,1]\}\]
is a nontrivial interval of real numbers. Thus, since $D_k$ is dense in $\R$, it follows that there is some $M_{n}\in[0,1]$ such that $f_{n-1}(w_k)+M_{n}g_{n}(w_k)\in D_k$. We define 
\[f_{n}(z)=f_{n-1}(z)+M_{n}g_{n}(z).\]

Let us show that (I)$_n$-(VII)$_n$ hold. It is trivial to see that (I)$_n$ and (II)$_n$ are true. For (III)$_n$, notice that because $M_n\in[0,1]$, and since (iii) and (III)$_{n-1}$ both hold, we have for all $x\in \R$,
\[f_n'(x)=f_{n-1}'(x)+M_{n}g'_{n}(x)\geq\frac{1}{2}+\frac{1}{2^{n-1}}-\frac{1}{2^{n}}=\frac{1}{2}+\frac{1}{2^{n}},\]
and thus $f_n:\R\to\R$ is a bijection. For (IV)$_n$ we have for all $z\in \C$,
\[|f_n(z)-f_{n-1}(z)|=M_{n}|g_{n}(z)|\leq \frac{1}{2^{n}}p(|z|),\]
where the last inequality follows since $M_n\in[0,1]$ and (ii) holds. Let us verify that (V)$_n$ holds.  From the definition of $f_n=f_{2k+1}$ and the way we chose $M_{n}$, it follows that $f_n(w_k)\in D_k$ (notice that $w_k \neq x_P$ by our case assumption). Thus (V)$_n$ holds. (VI)$_n$ holds trivially since $n$ is odd. To see that (VII)$_n$ holds, note that since $g_{n}(z)=0$ if $z\in \{x_P\}\cup A_{n-1}\cup f_{n-1}^{-1}(B_{n-1})$, it follows directly from the definition of $f_n$ that $f_n\restrict A_{n-1}=f_n\restrict A_{n-1}$ and $f_n^{-1}\restrict B_n=f_{n-1}^{-1}\restrict B_{n-1}$.

\textbf{Case \ref{sec_n_odd}.B:  $\boldsymbol{w_k\in \{x_P\}\cup A_{n-1}\cup f_{n-1}^{-1}(B_{n-1})}$}. Then we let $f_n=f_{n-1}$, $A_n=A_{n-1}\cup\{w_k\}$ and $B_n=B_{n-1}$.  Let us argue that this definition of $f_n$ satisfies (V)$_n$; the rest of (I)$_n$-(VII)$_n$ are easily seen to hold by the inductive hypothesis. Suppose $w_k\neq x_P$. Since the enumeration of $W$ is one-to-one we have $w_k\neq w_j$ for all $j\leq k-1$. Thus, for some $j\leq k-1$ we have $w_k=f_{n-1}^{-1}(w_j)$, and because $f_{n-1}(x_P)=y_P$, $f_{n-1}$ is injective and $w_k\neq x_P$, it follows that $w_j\neq y_P$. Since $2j+2\leq n-1$ and since it follows by our inductive assumptions (VII)$_{\ell}$ for $\ell\leq n-1$, that $f_{n-1}\restrict A_{2j+2}=f_{2j+2}\restrict A_{2j+2}$, we see that $w_k=f_{n-1}^{-1}(w_j)=f_{2j+2}^{-1}(w_j)\in D_j$. Then $D_j = D_k$ by (\ref{eq_KeyAssumption}) from page \pageref{eq_KeyAssumption}.  So, $f_n(w_k)=f_{n-1}(w_k)=w_j\in D_j=D_k$, and hence (V)$_n$ holds.

\subsection{When $n$ is even}\label{sec_n_even}  Now suppose $n=2k+2$ is even where $k>0$ and that $f_i$, $A_i$ and $B_i$ satisfying (I)$_i$-(VI)$_i$ have already been defined for $i\leq 2k+1$. We have
\[A_{2k+1}=\{w_0,\ldots,w_k\}\]
and
\[B_{2k+1}=\{w_0,\ldots,w_{k-1}\}.\]
We will define $f_n$, $A_n$ and $B_n$ in two cases as follows.

\textbf{Case \ref{sec_n_even}.A:  $\boldsymbol{f_{n-1}^{-1}(w_k)\notin \{x_P\}\cup A_{n-1}\cup f_{n-1}^{-1}(B_{n-1})}$}. Then we let $g_{n}$ be an entire function such that
\begin{enumerate}
\item[(i)] $(\forall z\in \C)$ $g_{n}(z)=0$ $\iff$ $z\in \{x_P\}\cup A_{n-1}\cup f_{n-1}^{-1}(B_{n-1})$,
\item[(ii)] $(\forall z\in \C)$ $|g_{n}(z)|\leq\frac{1}{2^{n}}p(|z|)$ and
\item[(iii)] $(\forall x\in \R)$ $g_{n}'(x)\geq -\frac{1}{2^{n}}$.
\end{enumerate}
For example, as in the case above where $n$ was odd, we could take
\[g_n(z)=\alpha_n(z-x_P)^{\beta_n}(z-w_0)\cdots(z-w_k)(z-f_{n-1}^{-1}(w_0))\cdots(z-f_{n-1}^{-1}(w_{k-1}))\]
satisfying (i)-(iii) by choosing $\alpha_n$ small enough and $\beta_n\in\{1,2\}$ so that the degree of $g_n$ is odd. By our inductive assumption about $f_{n-1}$ and by (iii) it follows that for any $M\in[0,1]$ and any $x\in\R$ we have \[f_{n-1}'(x)+Mg_{n}'(x)\geq\frac{1}{2}+\frac{1}{2^{n-1}}-\frac{1}{2^{n}}=\frac{1}{2}+\frac{1}{2^{n}}>0.\]
Thus the function $f_{n-1}+Mg_{n}:\R\to\R$ is a bijection. Let us argue that the set
\[\{(f_{n-1}+Mg_{n})^{-1}(w_k)\st M\in[0,1]\}\]
is a nontrivial interval of real numbers. It will suffice to show that $(f_{n-1}+g_n)^{-1}(w_k)\neq f^{-1}_{n-1}(w_k)$. Suppose $(f_{n-1}+g_n)^{-1}(w_k)= f_{n-1}^{-1}(w_k)$, then $f_{n-1}(f_{n-1}^{-1}(w_k))=w_k$ and $(f_{n-1}+g_n)(f^{-1}_{n-1}(w_k))=w_k$. This implies that the functions $f_{n-1}$ and $f_{n-1}+g_n$ are equal at the point $f_{n-1}^{-1}(w_k)$, and hence $g_n(f_{n-1}^{-1}(w_k))=0$, which contradicts (i) by our case assumption that $f_{n-1}^{-1}(w_k) \notin \{x_P\}\cup A_{n-1}\cup f_{n-1}^{-1}(B_{n-1})$.

Thus, since $D_k$ is dense in $\R$, it follows that there is some $M_{n}\in[0,1]$ such that $(f_{n-1}+M_{n}g_{n})^{-1}(w_k)\in D_k$. We fix such an $M_{n}$ and define 
\[f_{n}(z)=f_{n-1}(z)+M_{n}g_{n}(z).\]
We also let $A_n=A_{n-1}$
and $B_n=B_{n-1}\cup\{w_k\}$. The verification that (I)$_n$-(VII)$_n$ hold is straightforward and similar to the above; it is therefore left to the reader.

\textbf{Case \ref{sec_n_even}.B:   $\boldsymbol{f_{n-1}^{-1}(w_k)\in \{x_P\}\cup A_{n-1}\cup f_{n-1}^{-1}(B_{n-1})}$}, or equivalently $w_k\in \{y_P\}\cup f_{n-1}(A_{n-1})\cup B_{n-1}$. Then we define $f_n=f_{n-1}$. As in the odd case above, this definition of $f_n$ is easily seen to satisfy (I)$_n$-(V)$_n$ and (VII)$_n$. Let us check (VI)$_n$. Suppose $w_k\neq y_P$. Since the enumeration of $W$ is one-to-one we have $w_k\neq w_j$ for all $j\leq k-1$ and hence $w_k=f_{n-1}(w_j)$ for some $j\leq k$ where $w_j\neq x_P$. Since $2j+1\leq n-1$, it follows by our inductive assumptions (V)$_{\ell}$ for $\ell\leq n-1$, that $f_{n-1}\restrict A_{2j+1}=f_{2j+1}\restrict A_{2j+1}$ and $w_k=f_{n-1}(w_j)=f_{2j+1}(w_j)\in D_j$.  Then by \eqref{eq_KeyAssumption} from page \pageref{eq_KeyAssumption}, $D_j = D_k$.  So $f_n^{-1}(w_k) = f_{n-1}^{-1}(w_k) = w_j \in D_j = D_k$.

This concludes the proof of Theorem \ref{thm_StrengthenFranklin}.


\section{Proof of Theorem \ref{thm_StrengthenErdos}}\label{sec_Erdos}

To prove the $\Leftarrow$ direction of Theorem \ref{thm_StrengthenErdos}, assume that $\big\{ f_P \ : \ P \in \mathbb{R}^2 \big\}$ is a sparse analytic system, and consider the subcollection $\big\{ f_{(0,y)} \ : \ y \in \mathbb{R} \big\}$.  Since $f_{(0,y)}$ passes through the point $(0,y)$ and each $f_{(0,y)}$ is analytic, and hence continuous, it follows that for $y \ne y'$, $f_{(0,y)} \restriction (-\infty,0) \ne f_{(0,y')} \restriction (-\infty,0)$.  So
\[
\mathcal{F}:= \Big\{ f_{(0,y)} \restriction (-\infty,0) \ : \ y \in \mathbb{R} \Big\}
\]
is a continuum-sized collection of analytic functions on the common domain $D:=(-\infty,0)$.  Furthermore, given any $z \in D$, since $z \ne 0$ and the $f_P$'s formed a sparse analytic system, it follows that
\[
\{ f_{(0,y)} (z) \ : \ y \in \mathbb{R} \}
\]
is countable.  So $\mathcal{F}$ is a collection of analytic functions as in clause \ref{item_WetzelSystem} of Erd\H{o}s' Theorem \ref{thm_Erdos}.  So by that theorem, CH must hold.

To prove the $\Rightarrow$ direction of Theorem \ref{thm_StrengthenErdos}---which is heavily inspired by Erd\H{o}s' proof of Theorem \ref{thm_Erdos}---assume CH, and fix an enumeration $\langle w_\alpha \ : \ \alpha < \omega_1 \rangle$ of $\mathbb{R}$.  Fix any partition $\mathcal{D}$ of the reals into countable dense subsets of $\mathbb{R}$.\footnote{For example, define an equivalence relation $\sim$ on $\mathbb{R}$ by:  $x \sim y$ iff $y=rx$ for some nonzero $r \in \mathbb{Q}$.  Then the set of equivalence classes constitutes a partition of $\mathbb{R}$ into countable dense subsets of $\mathbb{R}$.  We thank Alex Misiats for pointing out this example (since our original draft used CH to get such a partition).}   For each $\alpha < \omega_1$, let $D_\alpha$ be the unique member of $\mathcal{D}$ containing $w_\alpha$.  Also fix an $\omega_1$-enumeration $\langle P_\alpha =(a_\alpha,b_\alpha) \ : \ \alpha < \omega_1 \rangle$ of $\mathbb{R}^2$.  

Fix an $\alpha < \omega_1$.  By Theorem \ref{thm_StrengthenFranklin}, there exists an entire $f_\alpha: \mathbb{C} \to \mathbb{C}$ such that:
\begin{enumerate}
    \item  $f_\alpha \restriction \mathbb{R}$ is a real analytic bijection with strictly positive derivative;
    \item $f_\alpha(a_\alpha) = b_\alpha$ (i.e., $f_\alpha \restriction \mathbb{R}$ passes through the point $P_\alpha$);
    \item For each $w_\xi$ in the countable set $W_\alpha:=\{ w_\xi \ : \ \xi < \alpha \}$:
    \begin{enumerate}
        \item\label{item_NE_forward}  if $w_\xi \ne a_\alpha$ then $f_\alpha(w_\xi) \in D_\xi$; and
        \item\label{item_NE_backward} if $w_\xi \ne b_\alpha$ then $f_\alpha^{-1}(w_\xi) \in D_\xi$.
    \end{enumerate}
\end{enumerate}

We claim that $\{ f_\alpha \restriction \mathbb{R} \ : \ \alpha < \omega_1 \}$ is a sparse analytic system, and the only nontrivial requirement to verify is that if $w \in \mathbb{R}$ then both
\[
A_w:=\big\{ f_\alpha(w) \ : \ \alpha < \omega_1 \text{ and }  w \ne a_\alpha  \big\}
\]
and
\[
B_w:=\big\{ f_\alpha^{-1}(w) \ : \ \alpha < \omega_1 \text{ and } w \ne b_\alpha \big\}
\]
are countable.  Say $w = w_\xi$; then 
\begin{align*}
A_w = A_{w_\xi} &  \subseteq  \underbrace{\{ f_\alpha(w_\xi) \ : \ \xi < \alpha < \omega_1 \text{ and } w_\xi \ne a_\alpha  \}}_{\subseteq D_\xi \text{, by \ref{item_NE_forward}}}   \cup  \underbrace{\{ f_\alpha(w_\xi) \ : \ \alpha \le \xi   \}}_{\text{countable because } \xi < \omega_1}   \end{align*}
and hence $A_w$ is countable.  Similarly,
\begin{align*}
B_w =B_{w_\xi} &  \subseteq  \underbrace{\{ f^{-1}_\alpha(w_\xi) \ : \ \xi < \alpha < \omega_1 \text{ and } w_\xi \ne b_\alpha  \}}_{\subseteq D_\xi\text{, by \ref{item_NE_backward}}}   \cup  \underbrace{\{ f^{-1}_\alpha(w_\xi) \ : \ \alpha \le \xi   \}}_{\text{countable because } \xi < \omega_1}
\end{align*}
and hence $B_w$ is countable.

\section{Proof of Theorem \ref{thm_AnswerBV_CH}}\label{sec_BV_solution}

The next lemma is the key connection between sparse analytic systems and predictors:  
\begin{lemma}\label{lem_SparseEquivRel}
Suppose $\mathcal{F}=\langle f_P \ : \ P \in \mathbb{R}^2 \rangle$ is a sparse analytic system.  Let $\sim$ be the equivalence relation on $\mathbb{R}$ generated by the set
\[
X:=\Big\{ (u,v) \in \mathbb{R}^2 \ : \ \exists P \in \mathbb{R}^2 \ \big( u \ne x_P \ \wedge \ v \ne y_P \ \wedge \ f_P(u)=v  \big) \Big\}.
\]
Then:
\begin{enumerate}
    \item\label{item_CtbleEquivClass} Each $\sim$-equivalence class is countable.
    \item\label{item_ClosedUnderEquiv} For every $P=(x_P,y_P) \in \mathbb{R}^2$ and every $z \in \mathbb{R}$:  if $z \ne x_P$ then $z \sim f_P(z)$.
\end{enumerate}

\end{lemma}

Before proving Lemma \ref{lem_SparseEquivRel}, we say how the proof of Theorem \ref{thm_AnswerBV_CH} is finished:  assuming CH, Theorem \ref{thm_StrengthenErdos} yields the existence of a sparse analytic system.  Let $\sim$ be the equivalence relation on $\mathbb{R}$ induced by the sparse analytic system via Lemma \ref{lem_SparseEquivRel}.  The properties of $\sim$ listed in the conclusion of Lemma \ref{lem_SparseEquivRel} satisfy the assumptions of Lemma 20 of Cox-Elpers~\cite{Cox_Elpers}; and that lemma tells us that if $S:=\mathbb{R}/\sim$ and 
\[
\mathcal{P}: {}^{\underset{\smile}{\mathbb{R}}} S \to S
\]
is any analytic-anonymous $S$-predictor,\footnote{Recall these notions were defined in Section \ref{sec_Intro}.} then $\mathcal{P}$ fails to predict the function $x \mapsto [x]_\sim$ for almost every $x \in \mathbb{R}$.\footnote{Strictly speaking, the \emph{statement} of \cite[Lemma 20]{Cox_Elpers} only implies that an analytic-anonymous predictor fails to predict $x \mapsto [x]_\sim$ on a \emph{positive}-measure set.  This is good enough to answer Question \ref{q_BV}, since such a predictor would not be good.  But the \emph{proof} of \cite[Lemma 20]{Cox_Elpers}---which was due essentially to Bajpai-Velleman~\cite{MR3552748}---shows that an analytic-anonymous predictor can successfully predict $x \mapsto [x]_\sim$ only for those $x$ lying in some fixed equivalence class, which in the context of Lemma \ref{lem_SparseEquivRel} is countable.  So analytic-anonymous predictors fail to predict $x \mapsto [x]_\sim$ almost everywhere in this situation.}  In particular, there is no good analytic-anonymous $S$-predictor.

\begin{proof}
(of Lemma \ref{lem_SparseEquivRel}): Part \eqref{item_ClosedUnderEquiv} holds because, by the definition of sparse analytic system, $f_P$ is injective and $f_P(x_P)=y_P$.  So if $z \ne x_P$ then $f_P(z) \ne y_P$; so not only is $z \sim f_P(z)$, but the pair $\big(z,f_P(z)\big)$ is an element of $X$.

To prove part \eqref{item_CtbleEquivClass}, since $X$ generates $\sim$, it suffices to prove that for every $z \in \mathbb{R}$, both
\[
z^\uparrow:=\{  v \in \mathbb{R} \ : \ (z,v) \in X  \} = \{ v \in \mathbb{R} \ :  \ \exists P \in \mathbb{R}^2  \ ( z \ne x_P \wedge v \ne y_P \wedge f_P(z)=v ) \}
\]
and
\[
z_\downarrow:=\{   u \in \mathbb{R} \ : \ (u,z) \in X \} =\{ u \in \mathbb{R} \ :  \ \exists P \in \mathbb{R}^2  \ ( u \ne x_P \wedge z \ne y_P \wedge f_P(u)=z ) \}
\]
are countable.  But 
\[
z^\uparrow \subseteq \{ f_P(z) \ : \ z \ne x_P \}
\]
and
\[
z_\downarrow \subseteq \{ f_P^{-1}(z) \ : \ z \ne y_P \}
\]
which are both countable, by definition of sparse analytic system.
\end{proof}

\section{Concluding Remarks}\label{sec_Concluding}

The notion of a sparse analytic system obviously generalizes to a sparse $\Gamma$-system for any $\Gamma \subseteq \text{Homeo}^+(\mathbb{R})$, and Lemma \ref{lem_SparseEquivRel} easily generalizes to such systems. In fact, Section 4 of Bajpai-Velleman~\cite{MR3552748} and Section 5 of Cox-Elpers~\cite{Cox_Elpers} can both be viewed as constructions, in ZFC alone, of sparse $\Gamma$-systems (with $\Gamma=$ ``increasing $C^\infty$ bijections" in \cite{MR3552748} and $\Gamma=$ ``increasing smooth diffeomorphisms" in \cite{Cox_Elpers}).

We have shown that CH implies a negative answer to Bajpai-Velleman's Question \ref{q_BV}, but it is open whether ZFC alone implies a negative solution.

\end{document}